\documentclass[11pt,reqno]{amsart} 

\usepackage{amssymb,latexsym}
\usepackage{cite} 

\usepackage[height=190mm,width=130mm]{geometry} 

\theoremstyle{plain}
\newtheorem{theorem}{Theorem}

\newtheorem{corollary}{Corollary}

\theoremstyle{definition}
\newtheorem{definition}{Definition}

\theoremstyle{remark}
\newtheorem{remark}{Remark}

\numberwithin{equation}{section} 

\begin{document}
\title[Closed ideals and bounded $\Delta$-weak approximate identities]{Closed ideals with bounded $\Delta$-weak approximate identities in some certain Banach algebras} 

\author{J. Laali}
\address{\textbf{Javad Laali}\\Kharazmi University\\Faculty of Mathematical and Computer Science\\ Department of Mathematics\\
 50 Taleghani Avenue, 15618\\
  Tehran\\ Iran}

\email{laali@khu.ac.ir}


\author{M. Fozouni}

\address{\textbf{Mohammad Fozouni}\\Gonbad Kavous University\\ Faculty of Sciences and Engineering\\Department of Mathematics\\
Gonbad Kavous\\ Iran}


\email{fozouni@hotmail.com}

\begin{abstract}
It is shown that a locally compact group $G$ is  amenable  if and only if some certain  closed ideals of the Fig\`{a}-Talamanca-Herz algebra $A_{p}(G)$ admit   bounded $\Delta$-weak approximate identities.   Also,  similar results  are obtained for the function algebras  $\mathfrak{L}A(G)$ and $C_{0}^{w}(G)$.
\end{abstract}




\maketitle

\section{Introduction and preliminaries}
Let  $A$ be a Banach algebra, $\Delta(A)$ be the character space of $A$, that is, the space of all non-zero homomorphisms from $A$ into $\mathbb{C}$ and $A^{*}$ be the dual space of $A$ consisting of all bounded linear functionals on $A$.
Throughout the paper,  $A$ is a commutative and semi-simple Banach algebra, hence $\Delta(A)$ is non-empty.

Let $\{e_{\alpha}\}$ be a net in  Banach algebra $A$. The net $\{e_{\alpha}\}$ is called,
\begin{enumerate}
  \item an \emph{approximate identity} if for each $a\in A$, $\|ae_{\alpha}-a\|\rightarrow 0$,
  \item a \emph{weak approximate identity} if for each $a\in A$, $|f(ae_{\alpha})-f(a)|\rightarrow 0$ for all $f\in A^{*}$,
  \item a \emph{$\Delta$-weak approximate identity} if for each $\phi\in \Delta(A)$, $|\phi(e_{\alpha})-1|\rightarrow 0$.
\end{enumerate}
\begin{definition}\label{maindef}
Let $A$ be a Banach algebra. A bounded $\Delta$-weak approximate identity for subspace $E\subseteq A$ is a bounded net $\{e_{\alpha}\}$ in $E$ such that for each $a\in E$,
\begin{equation*}
\lim_{\alpha}|\phi(ae_{\alpha})-\phi(a)|=0\hspace{0.5cm}(\phi\in \Delta(A)).
\end{equation*}
\end{definition}
For simplicity of notation,  let b.$\Delta$-w.a.i stand for a bounded $\Delta$-weak approximate identity.

It was proved that every Banach algebra $A$ with  a bounded $\Delta$-weak approximate identity has a bounded approximate identity (b.a.i) and  conversely \cite[Proposition 33.2]{Doran}.

The notion of a $\Delta$-weak approximate identity introduced and studied in \cite{Jones} where an example of a Banach algebra which has a $\Delta$-weak approximate identity but does not have an approximate identity, was given.
Indeed, if $S=\mathbb{Q}^{+}$ is the semigroup of positive rationales under addition, it was shown that the semigroup algebra $l^{1}(S)$ has a b.$\Delta$-w.a.i, but it does not have any bounded or unbounded approximate identity.

As the second example to see the difference between bounded approximate  and bounded $\Delta$-weak approximate identities, let $\mathbb{R}$ be the additive real line group and $1<p\leq\infty$. Put $S_{p}(\mathbb{R})=L^{1}(\mathbb{R})\cap L^{p}(\mathbb{R})$ and define the following norm;
\begin{equation*}
\|f\|_{S_{p}}=\max\{\|f\|_{1},\|f\|_{p}\}\hspace{0.5cm}(f\in S_{p}(\mathbb{R})).
\end{equation*}
Using \cite[Theorem 2.1]{Inoue}, we can see that $S_{p}(\mathbb{R})$  has a b.$\Delta$-w.a.i, but it has no b.a.i. Because we know that  $S_{p}(\mathbb{R})$ is a Segal algebra and it is well-known that a Segal algebra $S$ in $L^{1}(\mathbb{R})$ has a b.a.i if and only if $S=L^{1}(\mathbb{R})$. But it is clear that $S_{p}(\mathbb{R})\neq L^{1}(\mathbb{R})$.

These type of approximate identities have some interesting applications, for example; see \cite{Kamali, TH, Forrest}. In the past decades, B. E. Forrest studied the relations between the amenability of a group $G$ and closed ideals of $A(G)$  and $A_{p}(G)$ with a b.a.i; see \cite{Forrest1, F, F2, FKLS}, and the relations between some properties of  $G$ and closed ideals of $A(G)$   with a b.$\Delta$-w.a.i; see \cite{Forrest}.

 In this paper, we try to improve some of the  theorems in \cite{Forrest1, F, FKLS, GhL, Kamali} with changing b.a.i by b.$\Delta$-w.a.i. As an application, we give the converse of \cite[Corollary 4.2]{FKLS} due to B.  Forrest, E. Kaniuth, A. T. Lau and N. Spronk.
\section{Main results}
Let $G$ be a locally compact group. For $1<p<\infty$, let $A_{p}(G)$ denote the subspace of $C_{0}(G)$ consisting of all functions of the form $u=\sum_{i=1}^{\infty}f_{i}\ast \widetilde{g_{i}}$ where $f_{i}\in L^{p}(G)$, $g_{i}\in L^{q}(G)$, $1/p+1/q=1$, $\sum_{i=1}^{\infty}\|f_{i}\|_{p}\|g_{i}\|_{q}<\infty$ and $\widetilde{g}(x)=\overline{g(x^{-1})}$ for all $x\in G$. With the pointwise operation and the following norm,
\begin{equation*}
\|u\|_{A_{p}(G)}=\inf\{\sum_{i=1}^{\infty}\|f_{i}\|_{p}\|g_{i}\|_{q}: u=\sum_{i=1}^{\infty}f_{i}\ast \widetilde{g_{i}}\},
\end{equation*}
$A_{p}(G)$  is a Banach algebra called the Fig\`{a}-Talamanca-Herz algebra.
It is clear that $\|u\|\leq \|u\|_{A_{p}(G)}$ where $\|u\|$ is the uniform norm of $u\in C_{0}(G)$.  By \cite[Theorem 3]{Herz}, we know that
 $$\Delta(A_{p}(G))=\{\phi_{x} : x\in G\}=G,$$
where $\phi_{x}$ is defined by $\phi_{x}(f)=f(x)$ for each $f\in A_{p}(G)$.

The dual of the Banach algebra $A_{p}(G)$ is the Banach space $PM_{p}(G)$ consisting of all limits of convolution operators associated to bounded measures. Indeed, $PM_{p}(G)$ is the $w^{*}$-closure of $\lambda_{p}(L^{1}(G))$ in $B(L^{p}(G))$ where $\lambda_{p}$ is the left regular representation of $G$ on $L^{p}(G)$; see \cite{Derighetti} for more details.

The group $G$ is said to be \emph{amenable} if, there exists an $m\in L^{\infty}(G)^{*}$ such that $m\geq 0$, $m(1)=1$ and $m(L_{x}f)=m(f)$ for each $x\in G$ and $f\in L^{\infty}(G)$ where $L_{x}f(y)=f(x^{-1}y)$.
\begin{theorem}(Leptin-Herz)
Let $G$ be a locally compact group and $1<p<\infty$. Then $A_{p}(G)$ has a b.a.i if and only if $G$ is  amenable
\end{theorem}
The proof of the above theorem in the case $p=2$ is due to Leptin \cite{Leptin} and in general is due to Herz \cite{Herz}.

Forrest and  Skantharajah in \cite{Forrest} showed that if $G$ is a discrete group, then $A_{2}(G)=A(G)$ has a b.$\Delta$-w.a.i if and only if $G$ is amenable. Kaniuth and \"{U}lger in \cite[Theorem 5.1]{KU}, for the first time in our knowledge, announced that $A(G)$ has a b.$\Delta$-w.a.i if and only if $G$ is an amenable group, but the same result holds  for the Fig\`{a}-Talamanca-Herz algebras as follows. The proof is similar to the Fourier algebra case, but we give it for the convenience of reader.
\begin{theorem}\label{Lem: L-H}
 Let $G$ be a locally compact group and $1<p<\infty$. Then $A_{p}(G)$ has a b.$\Delta$-w.a.i if and only if $G$ is amenable.
\end{theorem}
\begin{proof}
Let $\{e_{\alpha}\}$ be a b.$\Delta$-w.a.i for $A_{p}(G)$ and $e\in A_{p}(G)^{**}$ be a $w^{*}$-cluster point of $\{e_{\alpha}\}$.
So, for each $\phi\in \Delta(A_{p}(G))=G$, we have
$$e(\phi)=\lim_{\alpha}\phi(e_{\alpha})=1,$$

because $\{e_{\alpha}\}$ is a b.$\Delta$-w.a.i for $A_{p}(G)$.
Therefore, by \cite[Proposition 2.8]{Ulger} $G$ is weakly closed in $PM_{p}(G)=A_{p}(G)^{*}$. Now, by \cite[Corollary 2.8]{Chou} we conclude that $G$ is an amenable group.
\end{proof}
The following corollary immediately follows from the Leptin-Herz Theorem and Theorem \ref{Lem: L-H}.
\begin{corollary}\label{Cor: main}
Let $G$ be a locally compact group and $1<p<\infty$. Then $A_{p}(G)$ has a b.$\Delta$-w.a.i. if and only if it has a b.a.i.
\end{corollary}
The following theorem is a key tool in the sequel.
\begin{theorem}\label{T:Ap identity} Let $A$ be a Banach algebra, $I$ be a closed two-sided  ideal of $A$ which has a b.$\Delta$-w.a.i and the quotient Banach algebra $A/I$ has a  b.l.a.i. Then $A$ has a b.$\Delta$-w.a.i.
\end{theorem}
\begin{proof}
Let $\{e_{\alpha}\}$ be a b.$\Delta$-w.a.i for $I$ and $\{f_{\delta}+I\}$ be a b.l.a.i for $A/I$.  We can assume that $\{f_{\delta}\}$ is bounded. Indeed, since $\{f_{\delta}+I\}$ is bounded, there exists a positive integer $K$ with $\|f_{\delta}+I\|<K$ for each $\delta$. So,  there exists $y_{\delta}\in I$ such that $\|f_{\delta}+I\|<\|f_{\delta}+y_{\delta}\|<K$. Put $f^{'}_{\delta}=f_{\delta}+y_{\delta}$. Clearly, $\{f^{'}_{\delta}+I\}$ is a b.l.a.i for $A/I$ which $\{f^{'}_{\delta}\}$ is  bounded.

Now, consider the bounded net $\{e_{\alpha}+f_{\delta}-e_{\alpha}f_{\delta}\}_{(\alpha,\delta)}$. For each $\phi\in \Delta(A)$ we have
\begin{equation*}
\phi(e_{\alpha}+f_{\delta}-e_{\alpha}f_{\delta})=\phi(e_{\alpha})+\phi(f_{\delta})(1-\phi(e_{\alpha}))\xrightarrow{(\alpha,\delta)}1.
\end{equation*}
Therefore, $A$ has a b.$\Delta$-w.a.i.
\end{proof}
Let $G$ be a locally compact group, $E$ be a closed non-empty subset of $G$ and $1<p<\infty$. Define
\begin{equation*}
I_{p}(E)=\{u\in A_{p}(G) : u(x)=0 \text{ for all } x\in E\}\cdot
\end{equation*}
The following result  improves \cite[Theorem 3.9]{Forrest1}.
\begin{theorem}\label{Th: 1}
Let $G$ be a locally compact group. Then the following assertions are equivalent.
\begin{enumerate}
\item $G$ is an  amenable group.
\item $\ker(\phi)$ has a b.$\Delta$-w.a.i for each $\phi\in \Delta(A_{p}(G))$.
\item $I_{p}(H)$ has a b.$\Delta$-w.a.i for some closed amenable subgroup $H$ of $G$.
\end{enumerate}
\end{theorem}
\begin{proof}
$(1)\Rightarrow (2)$: Let $G$ be an amenable group. Then $A_{p}(G)$ has a b.a.i by Leptin-Herz's Theorem. Now, the result follows from \cite[Corollary 2.3]{Kaniuth}.

$(2)\Rightarrow (3)$: Just take $H=\{e\}$, because we know that $I_{p}(\{e\})=\ker(\phi_{e})$.

$(3)\Rightarrow (1)$: Suppose that $I_{p}(H)$ for a closed amenable subgroup $H$ of $G$ has a b.$\Delta$-w.a.i. By  \cite[Lemma 3.19]{Monfared} we know that $A_{p}(H)$ is isometrically isomorphic to $A_{p}(G)/I_{p}(H)$. But $A_{p}(H)$ has a b.a.i, since $H$ is an amenable group. Therefore,
$A_{p}(G)/I_{p}(H)$ also has a b.a.i. Now, the result follows from Theorems \ref{T:Ap identity} and \ref{Lem: L-H}.
\end{proof}
Forrest in \cite[Lemma 3.14]{F} improved  \cite[Theorem 3.9]{Forrest1}. Indeed he showed that $G$ is an amenable group if for some closed proper subgroup $H$ of $G$, $I_{2}(H)$ has a b.a.i. Also, using the operator space structure of $A(G)$, it was shown in  \cite[Theorem 1.5]{FKLS} that $G$ is an amenable group only if $I_{2}(H)$ for some  closed subgroup $H$ of $G$ has a b.a.i.

Now, we give the following result which improves  \cite[Corollary 1.6]{FKLS} and Theorem \ref{Th: 1}. Our proof is a mimic of \cite[Lemma 3.14]{F}.
\begin{theorem}
Let $G$ be a locally compact group and $1<p<\infty$. Then the following are equivalent.
\begin{enumerate}
\item $G$ is an amenable group.
\item $I_{p}(H)$ has a b.$\Delta$-w.a.i for some proper closed subgroup $H$ of $G$.
\end{enumerate}
\end{theorem}
\begin{proof} In view of \cite[Corollary 4.2]{FKLS}, only $(2)\Rightarrow (1)$ needs proof.

Let $H$ be a proper closed subgroup of $G$ such that $I_{p}(H)$ has a b.$\Delta$-w.a.i. We will show that $H$ is an amenable group and this completes the proof by Theorem \ref{Th: 1}.

Since $H$ is a proper subgroup, there exists $x\in G\setminus H$. On the other hand, the mapping $I_{p}(H)\rightarrow I_{p}(xH)$ defined by $u\rightarrow L_{x}u$ is an isometric isomorphism, because for each $t\in G$ and $f\in A_{p}(G)$,  we have, $$L_{t}f\in A_{p}(G), \ \|L_{t}f\|_{A_{p}(G)}=\|f\|_{A_{p}(G)}.$$

 Therefore, $I_{p}(xH)$  has a b.$\Delta$-w.a.i which we denote it by $(u_{\alpha})$. For each $\alpha$, let $v_{\alpha}$ be the restriction of $u_{\alpha}$ to $H$. Using \cite[Theorem 1a]{Herz}, we conclude that $(v_{\alpha})$ is a bounded net in $A_{p}(H)$.

Let $\nu \in A_{p}(H)\cap C_{c}(H)$ and $K=\mathrm{supp} \ \nu\subseteq H$. Then there exists a neighborhood $V$ of $K$ in $G$ such that $V\cap xH=\emptyset$, because $K\cap xH=\emptyset$ (otherwise we conclude that $x$ is in $H$) and $G$ is  completely regular by \cite[Theorem 8.4]{Hew} and hence it is a regular topological space. Indeed, for each $y\in xH$, let $V_{y}$ be a neighborhood of $K$ such that $y\notin V_{y}$. So, $V=\cap_{y\in xH}V_{y}$ satisfies $V\cap xH=\emptyset$.

By \cite[Proposition 1, pp.34]{Derighetti} there is $u\in َA_{p}(G)$ such that $u(x)=1$ for each $x\in K$ and  $\mathrm{supp }\ u\subseteq V$, and by \cite[Theorem 1b]{Herz}, there exists a $v\in A_{p}(G)$ such that $v_{|H}=\nu$. Now, put $w=vu$. Since $V\cap xH=\emptyset$ and  $\mathrm{supp }\ u\subseteq V$, we have $w\in I_{p}(xH)$, and since $u(x)=1$ for each $x\in K$ and $K=\mathrm{supp} \ \nu$, we have $w_{|H}=\nu$.

Now, for each  $x\in H$, we have
\begin{align*}
\lim_{\alpha}|\phi_{x}(v_{\alpha}\nu)-\phi_{x}(\nu)|&=\lim_{\alpha}|v_{\alpha}(x)\nu(x)-\nu(x)|\\
&=\lim_{\alpha}|u_{\alpha}(x)w(x)-w(x)|=0.
\end{align*}
Therefore, $(v_{\alpha})$ is a b.$\Delta$-w.a.i for $A_{p}(H)$, since by \cite[Corollary 7, pp. 38]{Derighetti}, $A_{p}(H)\cap C_{c}(H)$ is dense in $A_{p}(H)$. Hence, by Theorem \ref{Lem: L-H}, $H$ is amenable.
\end{proof}
As an application of the above theorem, we give the following corollary which is the converse of \cite[Corollary 4.2]{FKLS}.
\begin{corollary}
Let $G$ be a locally compact group, $1<p<\infty$ and $H$ be a proper closed subgroup of $G$. If $I_{p}(H)$ has a b.a.i, then $G$ is amenable.
\end{corollary}
Ghahramani and Lau in \cite{GhL} introduced and studied a new closed ideal of $A(G)$. Indeed, let $G$ be a locally compact group and put $$\mathfrak{L}A(G)=L^{1}(G)\cap A(G)$$ with the norm $$|\|f\||=\|f\|_{1}+\|f\|_{A(G)}.$$

Clearly $\mathfrak{L}A(G)$   with  pointwise multiplication is a commutative Banach algebra with $\Delta(\mathfrak{L}A(G))=G$ and it is called  the \emph{Lebesgue-Fourier algebra} of $G$.

It was shown that $\mathfrak{L}A(G)$ has a b.a.i if and only if $G$ is a compact group \cite[Proposition 2.6]{GhL}.
 Now, we give the following result concerning the b.$\Delta$-w.a. identities of this Banach algebra.
\begin{theorem}\label{Th: LF}
Let $G$ be a locally compact group.
\begin{enumerate}
\item for each $x\in G$, $\ker(\phi_{x})\subseteq \mathfrak{L}A(G)$ has a b.$\Delta$-w.a.i.
\item $\mathfrak{L}A(G)$ has a b.$\Delta$-w.a.i.
\item  $G$ is amenable.
\item $G$ is compact.

Then $(1)\Rightarrow (2)\Rightarrow (3)$ and $(4)\Rightarrow (1)$ hold.
\end{enumerate}
\end{theorem}
\begin{proof} $(1)\Rightarrow (2):$ Follows from Theorem \ref{T:Ap identity}.

$(2)\Rightarrow (3):$ Let $(u_{\alpha})$ be a b.$\Delta$-w.a.i for $\mathfrak{L}A(G)$, Then
$(u_{\alpha})$ is a bounded net in $A(G)$. Now, the result follows from Theorem \ref{Lem: L-H} and this fact that $\Delta(\mathfrak{L}A(G))=G=\Delta(A(G))$.

$(4)\Rightarrow (1):$ Let $G$ be a compact group. Then by \cite[Proposition 2.6]{GhL}, we know that $\mathfrak{L}A(G)=A(G)$. On the other hand, by \cite[Example 2.6]{Kaniuth} for each $x\in G$, $A(G)$ is $\phi_{x}$-amenable. Therefore, the result follows from \cite[Proposition 2.2]{Kaniuth}.
\end{proof}
We do not know whether the implication $(3)\Rightarrow (2)$ in Theorem \ref{Th: LF} remains true?
\begin{remark} In \cite{G}, Granirer gave a $(p,r)$-version of the Lebesgue-Fourier algebra. Indeed, let $1<p<\infty, 1\leq r\leq \infty$ and put $A_{p}^{r}(G)=A_{p}(G)\cap L^{r}(G)$ with the norm,
$$\|u\|_{A_{p}^{r}(G)}=\|u\|_{r}+\|u\|_{A_{p}(G)}.$$

It was shown that $A_{p}^{r}(G)$ with pointwise multiplication is a commutative semi-simple Banach algebra such that $\Delta(A_{p}^{r}(G))=G$ and for all $1\leq r\leq \infty$, $A_{p}^{r}(G)=A_{p}(G)$ if $G$ is a compact group; see \cite[Theorem 1, Theorem 2]{G}.

Therefore, in view of Theorems \ref{Lem: L-H} and \ref{T:Ap identity}, Theorem \ref{Th: LF} remains true if we replace $\mathfrak{L}A(G)$ with $A_{p}^{r}(G)$.
\end{remark}
\begin{remark} Runde in \cite{Runde}, by using  the canonical operator space structure of $L^{p}(G)$, introduced and studied  the algebra O$A_{p}(G)$ for $1<p<\infty$, the \emph{operator Fig\`{a}-Talamanca-Herz algebra}.  It was shown that $A_{p}(G)\subseteq \mathrm{O}A_{p}(G)$ \cite[Remark 4. pp. 159]{Runde}  and  O$A_{p}(G)$ has a b.a.i if and only if $G$ is an amenable group \cite[Theorem 4.10]{Runde}.
That would be an interesting question: Are the preceding results  remain true if we replace $A_{p}(G)$ by O$A_{p}(G)$?
\end{remark}

Now, let $G$ be a locally compact group
 and $w: G\rightarrow \mathbb{R}$ be an upper semicontinuous function such that $w(x)\geq 1$ for each $x\in G$. Put
\begin{equation*}
C_{0}^{w}(G)=\{f\in C(G) : fw\in C_{0}(G)\}.
\end{equation*}
It is clear that $C_{0}^{w}(G)$ with pointwise operation and weighted supremum norm defined by
\begin{equation*}
\|f\|_{w}=\sup_{x\in G} |f(x)|w(x)\quad (f\in C_{0}^{w}(G)),
\end{equation*}
is a commutative Banach algebra such that $\Delta(C_{0}^{w}(G))=\Delta(C_{0}(G))=G$; see \cite[Section 4.3]{Kamali}. 
\begin{theorem}\label{Th: w}
Let $G$ be a locally compact group and $t\in G$.  Then the following  are equivalent.
\begin{enumerate}
\item $\ker(\phi_{t})$ has a b.a.i.

\item $\ker(\phi_{t})$ has a b.$\Delta$-w.a.i.

\item $C_{0}^{w}(G)$ has a b.$\Delta$-w.a.i.

\item w is bounded.
\end{enumerate}
\end{theorem}
\begin{proof} $(1)\Rightarrow (2):$ This part is clear.  Applying Theorem \ref{T:Ap identity}, we conclude $(2)\Rightarrow (3)$. In view of \cite[Corollary 4.7]{Kamali}, we have $(3)\Rightarrow (4)$. Therefore, we only show $(4)\Rightarrow (1)$.

 Suppose that $w$ is bounded by $M>0$. Let  $\mathcal{V}=\{V_{\alpha}\}_{\alpha\in \Gamma}$ be a neighborhood base at $t$ directed by the reverse inclusion. For each $\alpha\in \Gamma$, by the Urysohn Lemma, there exists a continuous function $f_{\alpha}:X\rightarrow [0,1]$ such that $f_{\alpha}(t)=1$ and $\mathrm{supp}f_{\alpha}\subseteq V_{\alpha}$.

Let $\epsilon>0$ and $g\in C_{0}^{w}(G)$. For $\epsilon^{'}=\epsilon/2M$, there exists a neighborhood  $V$ of $t$ such that,
\begin{equation*}
|g(y)-g(t)|<\epsilon^{'}\quad(y\in V).
\end{equation*}
On the other hand, let $\alpha_{0}\in \Gamma$ be such that $V_{\alpha_{0}}\subseteq V$. Therefore, there exists $x_{0}\in G$ such that
\begin{align*}
\|gf_{\alpha}-\phi_{t}(g)f_{\alpha}\|_{w}&=\sup_{x\in G}|g(x)f_{\alpha_{0}}(x)-g(t)f_{\alpha_{0}}(x)|w(x)\\
 &<|g(x_{0})f_{\alpha_{0}}(x_{0})-g(t)f_{\alpha_{0}}(x_{0})|w(x_{0})+\epsilon/2\\
 &<M\epsilon^{'}+\epsilon/2=\epsilon.
\end{align*}
Therefore, for each $g\in C_{0}^{w}(G)$,  $\|gf_{\alpha}-g(t)f_{\alpha}\|_{w}\longrightarrow 0$.
 Hence, by  \cite[Theorem 1.4, Proposition 2.2]{Kaniuth} and \cite[Corollary 3.6]{AK} we conclude that $\ker(\phi_{t})$ has a b.a.i.
\end{proof}

\section*{Acknowledgement} 
The authors are grateful to the referee of the paper for his (her) invaluable comments which improved the presentation of the paper, especially for giving a shorter proof for Theorem \ref{T:Ap identity}.

\bibliography{mmnsample}
\bibliographystyle{mmn}

\end{document}